\documentclass[12pts]{amsart}

\usepackage{mathpazo}
\usepackage[utf8]{inputenc}
\usepackage{amsmath}
\usepackage{amsthm}
\usepackage{amssymb}
\usepackage{amscd}
\usepackage{amsfonts}

\usepackage[utf8]{inputenc}
\usepackage{hyperref}			
\usepackage{cleveref}			
\usepackage{dsfont}	
\usepackage{mathtools}
\usepackage{enumitem}
\usepackage{dblfloatfix}
\usepackage{subfig}
\usepackage{biblatex} 
\usepackage[centering,margin=2cm]{geometry}

\newtheorem{theo}{Theorem}

\newtheorem{prop}{Proposition}
\newtheorem{rem}{Remark}
\theoremstyle{definition}

\newcommand{\eps}{\varepsilon}
\newcommand{\R}{\mathds{R}}
\newcommand{\N}{\mathds{N}}

\newcommand{\sig}{\Sigma}

\begin{document}

\title{Sharp estimates for the number of limit cycles in discontinuous generalized Liénard equations}
\author[T. M. P. de Abreu, R. M. Martins]{Tiago M. P. de Abreu${}^1$ and Ricardo M. Martins${}^2$}
	
\address{Departamento de Matem\'{a}tica, Universidade
Estadual de Campinas, Rua S\'{e}rgio Baruque de Holanda, 651, Cidade Universit\'{a}ria Zeferino Vaz, 13083--859, Campinas, SP,Brazil} \email{TMPA: tmpabreu@ime.unicamp.br, RMM: rmiranda@unicamp.br}
	
\subjclass[2010]{34A36,34C28,34C15,37C27,37G15}
	
\keywords{piecewise smooth differential equations, averaging theory, Lienard systems}
	
\maketitle

\begin{abstract}
In this paper, we study the maximum number of limit cycles for the piecewise smooth system of differential equations $\dot{x}=y, \ \dot{y}=-x-\eps \cdot (f(x)\cdot y +\text{sgn}(y)\cdot g(x))$. Using the averaging method, we were able to generalize a previous result for Liénard systems. In our generalization, we consider $g$ as a polynomial of degree $m$. We conclude that for sufficiently small values of $|\epsilon|$, the number $\left[\frac{n}{2}\right]+\left[\frac{m}{2}\right]+1$ serves as a lower bound for the maximum number of limit cycles in this system, which bifurcates from the periodic orbits of the linear center $\dot{x}=y$, $\dot{y}=-x$. Furthermore, we demonstrate that it is indeed possible to achieve such a number of limit cycles.

\end{abstract}

\section{Introduction}

The existence of periodic orbits and limit cycles for differential systems is a very active topic of research, mainly due to its relation to the Hilbert 16th Problem, and the existence of such objects in systems of piecewise smooth differential equations has been studied recently by several authors. For cases in $\R^2$, there are many works that determine the maximum number of limit cycles for a given class of vector fields and different separation manifolds (see \cite{caldas2022}, \cite{llibre-zhang}, \cite{novaes-llibre-teixeira}). There are also studies of piecewise smooth differential systems in $\R^3$ (see \cite{Llibre-Teixeira-R3}, \cite{llibre-tonon-velter}), for which the discontinuity manifold is a plane.
\par Moreover, there are also many open cases. For instance, for piecewise linear systems with a straight line of discontinuity, \cite{huan-yang} and \cite{llibre-ponce} prove the existence of a case with three nested limit cycles, and \cite{llibre-li} provides a very detailed description of different configurations for such systems and their respective lower bounds for the maximum number of limit cycles. However, as far as we know, the maximum number of limit cycles in piecewise linear systems with a straight line of discontinuity is still an open question. Nevertheless, droping the assumption that the discontinuity is a straight line, \cite{braga-melo} and \cite{douglas-ponce} show us that there exists piecewise linear systems having exactly $n$ limit cycles, for any $n \in \N$.
\par The generalized Liénard differential equation
\begin{equation}\label{intro-lienard}
    \ddot{x}+f(x)\dot{x}+g(x)=0,\end{equation}
 introduced in \cite{lienard1928}, where $f(x)$ and $g(x)$ are polynomial functions on the variable $x$, is a subclass of polynomial vector fields with applications on modelling oscilattion circuits that have been extensively studied in the past years. 

 This differential system has been extensively studied in the literature \cite{ASHEGHI2015120,YANG2011269,YU2007617,YU20061048}. In 1976, Lins, de Melo and Pugh \cite{lins-melo-pugh} conjectured that for $g(x)=x$ and $f(x)$ of degree $n\geq 1$, the system \ref{intro-lienard} has at most $\left[\frac{n-1}{2}\right]$ limit cycles --- where the brackets denote the integer part function. This was proved to be true for $n=4$ \cite{llibre-li}. For $n\geq 6$, Maesschalck and Dumortier proved in 2011 \cite{MAESSCHALCK} that such equations can have $\left[\frac{n-1}{2}\right]+2$ limit cycles, contradicting the conjecture from \cite{lins-melo-pugh}. The conjecture for $n=5$ remains open.
\par  One of the works that started the study of limit cycles on discontinuous Liénard polynomial differential systems was done by Llibre and Teixeira and is summarized by the paper \cite{llibre-teixeira}. In their work, they proposed the study of the system
\begin{align}\label{llibre-teixeira-syst}
    \begin{matrix*}[l]
        \dot{x}=y+\eps\cdot \text{sgn}(g_m(x,y))\cdot f(x), \\
        \dot{y}=-x,
    \end{matrix*}
\end{align}
where $f(x)$ is a polynomial of degree $n$ and the zero set of the function $\text{sgn}(g_m(x,y))$, $m\in\{0,2,4,6...\}$, is the union of $m/2$ distinct straight lines passing through the origin, dividing the plane in sectors of angles $2\pi/m$. They managed to prove, using the periodic averaging method for regularized discontinuous systems, that for $m=0,2$ and 4 the lower bounds for the maximum number of limit cycles of system \eqref{llibre-teixeira-syst} are, respectively, $\left[\frac{n-1}{2}\right]$, $\left[\frac{n}{2}\right]$ and $\left[\frac{n-1}{2}\right]$. They also left unproved the conjecture that, for $m\geq 6$, a lower bound for the maximum number of limit cycles of this system should be $\left[\frac{1}{2}\left(n-\frac{m-2}{2}\right)\right]$. This was then proven to be true in \cite{dong-liu}.
\par In \cite{miranda}, the authors studied the number of limit cycles in the system
\begin{align}\label{miranda-syst}
    \begin{matrix*}[l]
        \dot{x}=y, \\
        \dot{y}=-x-\eps(f(x)\cdot y+ \text{sgn}(y)(k_1x+k_2)),
    \end{matrix*}
\end{align}
where $f$ is a polynomial of degree $n\in \N$ and $k_1,k_2 \in \R$. The approach chosen by them was also based on the ideas of the regularization method --- specifically by introducing a piecewise linear function of the type
$$\varphi_w(y)=\begin{cases}-1 &\text{ if } y<-w,\\ \frac{y}{w} &\text{ if } -w<y<w,\\ 1 &\text{ if } y>w.\end{cases}.$$
\par The original system \eqref{miranda-syst} is then replaced by
\begin{align}\label{miranda-syst2}
    \begin{matrix*}[l]
        \dot{x}=y, \\
        \dot{y}=-x-\eps(f(x)\cdot y+ \varphi_w(y)(k_1x+k_2)),
    \end{matrix*}
\end{align}
and it is easy to see that, taking $w\rightarrow 0$, $\varphi_w(y)\rightarrow \text{sgn}(y)$. The main result of \cite{miranda} is the following theorem:
\begin{theo}\label{martins-mereu}\emph{(Martins-Mereu)}
 For every $n \geq 1$ and $|\eps|$ sufficiently small, the maximum number of limit cycles of the system \eqref{miranda-syst2} bifurcating from the periodic orbits of the linear center $\dot{x}=y$, $\dot{y}=-x$ is $[n/2]+1$. Moreover, there are systems \eqref{miranda-syst2} having exactly $[n/2]+1$ limit cycles.
\end{theo}
\par By proving this case, the similar problem for the differential system \eqref{miranda-syst} comes as a corollary, taking the limit $w \rightarrow 0$. In this paper, we generalize this theorem replacing $k_1x+k_2$ by an arbitrary real polynomial $g(x)$ of degree $m\geq 1$. Our proof will not pass through the regularization process --- in one hand, this should ease the process of calculating the averaged function but, in the other hand, we need to take a more careful look at the conditions of the functions we are averaging. Our main result is the following:
\begin{theo}
\emph{(Main Result)} Let $f(x)$ and $g(x)$ be real polynomials of degrees $n \geq 1$ and $m \geq 1$, respectively, and consider the system
  \begin{align}\label{main-syst}
        \begin{matrix*}[l]\Dot{x}=y,\\ \Dot{y}=-x-\eps \cdot (f(x)\cdot y +\text{sgn}(y)\cdot g(x)).\end{matrix*}
    \end{align}
Then, for $|\eps|$ sufficiently small, the number $\left[\frac{n}{2}\right]+\left[\frac{m}{2}\right]+1$ is a lower bound to the maximum number of limit cycles of the system \eqref{main-syst} bifurcating from the periodic orbits of the linear center $\dot{x}=y$, $\dot{y}=-x$. Moreover, we can choose $f$ and $g$ such that this number of cycles is indeed achieved.
\end{theo}
\begin{rem}
 Even though we aim to find a \emph{maximum number} of limit cycles of a system, we can only guarantee a \emph{lower bound of the maximum number} of limit cycles, in the sense of it being the maximum number of cycles \emph{that can be found through first-order averaging}.
\end{rem}

\section{Preliminaries}
\par Let $U$ be an open neighbourhood of $0$ where the vector field given by a piecewise smooth differential system is defined, and let $\sig$ be the hypersurface along which the discontinuities occur. Since any embedded hypersurface is locally the inverse image of a regular value, let $\sig=f^{-1}(0)\cap U$, where $f$ is the \emph{germ} of a $C^r$ function with $r>1$ and which has 0 as a regular value. Here the \emph{germ of a function} refers to the equivalence class of all functions which are locally equal to one another, see \cite{guardia-seara-teixeira} for more details.
\par Note that the hypersurface $\sig$ splits $U$ into the following open sets:
\begin{align*}
    \sig^+=\{(x,y)\in U: f(x,y)>0\} \text{ and }  \sig^-=\{(x,y)\in U: f(x,y)<0\}.
\end{align*}
\par We may then define the germs of a discontinuous vector field as
\begin{align}\label{filippovvect}
    Z(x,y)= \begin{cases}
    X(x,y), \text{ if } (x,y)\in \sig^+ \\
    Y(x,y), \text{ if } (x,y)\in \sig^-
    \end{cases},
\end{align}
and we denote the vector field by $Z=(X,Y)$.
\par The trajectories on $\sig^+$ and $\sig^-$ that do not reach $\sig$ can be defined by the vector fields $X$ and $Y$ in the same way as for smooth systems. Troubles may arise from what happens along $\sig$, hence we must take a careful look into these possibilities.
\par Let $Xf(p)=X(p)\cdot \nabla f(p)$ denote the Lie derivative of $f$ with respect to $X$ at the point $p$, and define:
\begin{align*}
    & \sig^c=\{p\in \sig: Xf(p)\cdot Yf(p)>0\},\\
    & \sig^s=\{p\in \sig: Xf(p)<0, Yf(p)>0\},\\
    & \sig^e=\{p\in \sig: Xf(p)>0, Yf(p)<0\},
\end{align*}
which we name \emph{crossing region, sliding region} and \emph{escaping region}, respectively. 

If $Xf(p)=0$ or $Yf(p)=0$, then $p$ is a \emph{tangency point}; we assume that these tangency points are isolated in $\sig$. 
\par Note that if $p\in \sig^c$, i.e. p is in a crossing region, then we can define the trajectory through $p$ by simply matching the trajectories defined by $X$ and $Y$. However, if $p \in \sig^s$ or $\sig^e$, the trajectory can't be defined so directly; in this case, one may use the Filippov convention, but this case will not be studied in this paper.

\par Now, the differential system that we're interested on averaging shall have the following form
  \begin{align}\label{ns-aver-syst}
      \dot{x}(t)=\eps F(t,x)+\eps^2 R(t,x,\eps),
  \end{align}
  with
  \begin{align*}
   & F(t,x)= \begin{cases}
    F^1(t,x), \text{ if } x\in S_1, \\
    F^2(t,x), \text{ if } x\in S_2,
    \end{cases}\\
     & R(t,x,\eps)= \begin{cases}
    R^1(t,x,\eps), \text{ if } x\in S_1, \\
    R^2(t,x,\eps), \text{ if } x\in S_2,
    \end{cases}
\end{align*}
where $F^{1,2}:I\times S_{1,2}\rightarrow \R^n$, $R^{1,2}:I\times S_{1,2} \times (-\eps_0,\eps_0)\rightarrow \R^n$ and $S_{1,2}=D\cap\sig^{\pm}$. Or, alternatively, using the function $\mathcal{X}(A)=\begin{cases}
    1, \text{ if } x\in A \\
    0, \text{ if } x\notin A
    \end{cases}$,
    we can rewrite those as:
  \begin{align*}
      & F(t,x)=\mathcal{X}_{S_1}F^1(t,x)+\mathcal{X}_{S_2}F^2(t,x),\\
        & R(t,x,\eps)=\mathcal{X}_{S_1}R^1(t,x,\eps)+\mathcal{X}_{S_2}R^2(t,x,\eps).
  \end{align*}
  \par Then, if $F^{1,2}$ and $R^{1,2}$ are $T$-periodic functions, we can define the averaged function $F_0(z)$ as 
  \begin{align}\label{ns-averaged-func}
   F_0(z)=\frac{1}{T}\int_0^T F(t,z)dt,
  \end{align}

  \par In \cite{llibre-mereu-novaes}, the authors proved the following result:
  \begin{theo}\label{aver-theo-2}
   \emph{(First order averaging for discontinuous systems, \cite{llibre-mereu-novaes})} Assuming the following hypothesis:
  \begin{enumerate}[label=H\arabic*]
  \item There exists an open bounded set $C\subset D$ such that, for each $z\in \overline{C}$, the curve $\{(t,z):t\in I= \mathds{S}^1=\R/T \}$ reaches transversely the set $\sig$ and only at generic points of discontinuity;
      \item For $j=1,2$, the continuous functions $F^j$ and $R^j$ are locally Lipschitz with respect to $x$, and T-periodic with respect to the time variable $t$;
      \item For $a\in C$ with $F_0(a)=0$, there exists a neighborhood $U\subset C$ of $a$ such that $F_0(z)\neq 0$ for all $z \in \overline{U}\setminus \{a\}$ and $d_B(F_0,U,0)\neq 0$. 
  \end{enumerate}
 Then for $|\eps|\neq 0$ sufficiently small, there exists a T-periodic solution $x(t,\eps)$ of system \eqref{ns-aver-syst} such that $x(0,\eps)\rightarrow a$ as $\eps \rightarrow 0$.
  \end{theo}

\par Theorem \ref{aver-theo-2} stabilishes the correspondence between the zeros of the averaged function and the closed orbits of the system \ref{ns-aver-syst}, so that roughly we can replace the problem of finding limit cycles by the problem of finding zeros of a polynomial. In order to deal with these polynomials, we introduce the following theorem, as stated in \cite{llibre-teixeira}:
\begin{theo}\label{theodescartes}\emph{(Descartes' Theorem)}
    Consider the real polynomial $p(x)=a_{i_1}x^{i_1}+...+a_{i_r}x^{i_r}$, with $r>1$, $0 \leq i_1<...<i_r$ and the numbers $a_{i_j}$ are not simultaneously zeros for $j\in\{1,2,...,r\}$. If $a_{i_j}\cdot a_{i_{j+1}}<0$, we say that they have a variation of sign. If the number of variations of signs is $m$, then $p(x)$ has at most $m$ positive real roots. Moreover, it's always possible to choose the coefficients of $p(x)$ in such a way that $p(x)$ has exactly $r-1$ positive real roots.
   \end{theo}
   
   \par With these tools, we may now prove the main theorem of this work.
   
\section{Proof of the Main Result}
\par Let $f(x)=\sum_{i=0}^{n}a_ix^i$ and $g(x)=\sum_{j=0}^m b_jx^j$. We start rewriting the system in polar coordinates
   \begin{align*}
        \begin{matrix}
            \dot{x} = \dot{r}\cos{\theta}-r\sin{\theta}\cdot\dot{\theta} \\
            \dot{y}=\dot{r}\sin{\theta}+r\cos{\theta}\cdot\dot{\theta} 
        \end{matrix} =\begin{pmatrix}\cos{\theta}&-r\sin{\theta}\\ \sin{\theta} & r\cos{\theta}\end{pmatrix}\cdot \begin{pmatrix}\dot{r}\\ \dot{\theta}\end{pmatrix},
    \end{align*}
then inverting the coordinate change matrix:
    \begin{align*}
        \begin{pmatrix}\dot{r}\\ \dot{\theta}\end{pmatrix}= &\frac{1}{r} \begin{pmatrix}r\cos{\theta}&r\sin{\theta}\\- \sin{\theta} & \cos{\theta}\end{pmatrix} \begin{pmatrix}
            r\sin{\theta}  \\ -r\cos{\theta} -\eps \cdot( f(r\cos{\theta})\cdot r \sin{\theta}+\text{sgn}(r\sin{\theta})\cdot g(r \cos{\theta}))
        \end{pmatrix}\\
            &=\begin{pmatrix}
            -\eps\cdot \sin{\theta}\cdot ( f(r\cos{\theta})\cdot r \sin{\theta}+\text{sgn}(r\sin{\theta})\cdot g(r \cos{\theta}))  \\ -1 -\frac{\eps}{r} \cdot( f(r\cos{\theta})\cdot r \sin{\theta}+\text{sgn}(r\sin{\theta})\cdot g(r \cos{\theta}))
        \end{pmatrix}.
    \end{align*}
\par To ease our calculations, put $a=( f(r\cos{\theta})\cdot r \sin{\theta}+\text{sgn}(r\sin{\theta})\cdot g(r \cos{\theta}))$. Then, admitting $\theta$ as the new independent variable, we have:
\begin{align*}
    \frac{d r}{d \theta}=\frac{-\eps \cdot \sin{\theta}\cdot a}{-1-\frac{\eps}{r}\cdot \cos{\theta}\cdot a}=\phi(\eps).
\end{align*}
\par Expanding the series of $\phi$ around $\eps=0$ we get
\begin{align*}
    \phi(\eps)=\phi(0)+\phi'(0)\cdot \eps +\mathcal{O}(\eps^2).
\end{align*}
\par Note that $\phi(0)=0$, and 
$$\phi'(0)=\left.\frac{-\sin{\theta}\cdot a(-1-\frac{\xi}{r}\cdot \cos{\theta} \cdot a)-(-\xi \sin{\theta}\cdot a \cdot \frac{\cos{\theta}\cdot a}{r})}{(-1-\frac{\xi}{r}\cdot \cos{\theta}\cdot a)^2}\right|_{\xi=0}=\sin{\theta}\cdot a ,$$
hence
$$\frac{d r}{d \theta}= \sin{\theta}\cdot ( f(r\cos{\theta})\cdot r \sin{\theta}+\text{sgn}(r\sin{\theta})\cdot g(r \cos{\theta}))\cdot \eps + \mathcal{O}(\eps^2).$$
\par Writing $F(r,\theta)=\sin{\theta}\cdot ( f(r\cos{\theta})\cdot r \sin{\theta}+\text{sgn}(r\sin{\theta})\cdot g(r \cos{\theta}))$ and $\mathcal{O}(\eps^2)=R(r,\theta)\cdot \eps^2$, define the averaged function $F_0(r)$ as the following integral:
\begin{align*}
    F_0(r)=\frac{1}{2\pi}\int_0^{2\pi}( f(r\cos{\theta})\cdot r \sin^2{\theta}+\text{sgn}(r\sin{\theta})\cdot \sin{\theta}\cdot g(r \cos{\theta})) d \theta .
\end{align*}
\par Note that $\text{sgn}(r\sin{\theta})=\text{sgn}(\sin{\theta})$, since $r>0$, which yields the value $+1$ when $\theta \in (0,\pi)$ and $-1$ when $\theta \in (\pi,2\pi)$, then:
\begin{align*}
    F_0(r)=\frac{1}{2\pi}\left[\int_0^{2\pi}( f(r\cos{\theta})\cdot r \sin^2{\theta} d\theta+\int_0^{\pi} \sin{\theta}\cdot g(r \cos{\theta})) d \theta - \int_{\pi}^{2\pi} \sin{\theta}\cdot g(r \cos{\theta})) d \theta \right],
\end{align*}
and since 
$$\int_{\pi}^{2\pi} \sin{\theta}\cdot g(r \cos{\theta})) d \theta = \int_{\pi}^{0} \sin{\theta}\cdot g(r \cos{\theta})) d \theta = -\int_{0}^{\pi} \sin{\theta}\cdot g(r \cos{\theta})) d \theta,  $$ the averaged function can be written as:
\begin{align*}
    F_0(r)=\frac{1}{2\pi}\left[\int_0^{2\pi}( f(r\cos{\theta})\cdot r \sin^2{\theta} d\theta +2\cdot \int_0^{\pi} \sin{\theta}\cdot g(r \cos{\theta})) d \theta\right].
\end{align*}
\par Let $I$ and $J$ denote the following integrals:
\begin{align*}
&I=\int_0^{2\pi}( f(r\cos{\theta})\cdot r \sin^2{\theta}, \\
&J=\int_0^{\pi} \sin{\theta}\cdot g(r \cos{\theta})) d \theta.
\end{align*}
\par For the first integral we have:
\begin{align*}
    I=\int_0^{2\pi} \sum_{i=0}^{n}a_i r^i\cos^i{\theta}\cdot r \sin^2{\theta} d\theta =\sum_{i=0}^{n}a_i r^{i+1} \cdot \int_0^{2\pi}\cos^i{\theta}\sin^2{\theta} d\theta .
\end{align*}
\par As in \cite{miranda}, we use the following formulas: 
\begin{align*}
    & \int_0^{2\pi} \cos^{2k+1}\theta \cdot \sin^2{\theta} d\theta =0, k=0,1,2... \\
     & \int_0^{2\pi} \cos^{2k}\theta \cdot \sin^2{\theta} d\theta =\pi \alpha_k\neq 0, k=0,1,2... 
\end{align*}
thus
\begin{align*}
    I=\sum_{i=0}^{\left[\frac{n}{2}\right]} \pi a_{2i} \alpha_{2i} r^{2i+1}. 
\end{align*}
\par Hence $I$ is a polynomial formed exclusively by odd exponents of $r$, and with $\left[\frac{n}{2}\right]+1$ terms. 
\par For the second integral, notice that
\begin{align*}
    J=\int_0^{\pi} \sum_{j=0}^{m}b_j r^j\cos^j{\theta}\cdot \sin{\theta} d\theta = \sum_{j=0}^{m}b_j r^j\int_0^{\pi} \cos^j{\theta}\cdot \sin{\theta} d\theta .
\end{align*}
\par In order to evaluate this integral, consider the change of variables $u=\cos{\theta}$, $du=\sin{\theta} d\theta$ then:
$$\int_0^{\pi} \cos^j{\theta}\cdot \sin{\theta} d\theta =\int_{1}^{-1}u^j du = \frac{(-1)^{j+1}-1}{j+1}.$$
\par If $j$ is odd, then $j+1$ is even and $\frac{(-1)^{j+1}-1}{j+1}=0$; but, when $j$ is even, i.e. when $j+1$ is odd, then $\frac{(-1)^{j+1}-1}{j+1}=-\frac{2}{j+1}\neq 0$. Hence
\begin{align*}
    J=\sum_{j=0}^{\left[\frac{m}{2}\right]}\frac{-2 b_{2j} r^{2j}}{2j+1},
\end{align*}
which means that $J$ is a polynomial with $\left[\frac{m}{2}\right]+1$ monomials, formed by even exponents of $r$. 
\par Since $ F_0(r)=\frac{1}{2\pi}(I+2J)$ and $I$ and $J$ don't have any powers of $r$ in common, it follows that $F_0(r)$ is a polynomial with $\left[\frac{n}{2}\right]+\left[\frac{m}{2}\right]+2$ terms. Therefore, by the Descartes' Theorem, $F_0(r)$ can have at most $\left[\frac{n}{2}\right]+\left[\frac{m}{2}\right]+1$ roots; then, if we could apply Theorem \ref{aver-theo-2}, it would follow that the system \eqref{main-syst} can have $\left[\frac{n}{2}\right]+\left[\frac{m}{2}\right]+1$ limit cycles bifurcating from the linear center. Moreover, since we can choose the coefficients from $F_0(r)$ such that it has exactly $\left[\frac{n}{2}\right]+\left[\frac{m}{2}\right]+1$ roots, this maximum number of limit cycles can indeed be achieved. Thus, to complete this proof, all we need to do is verify that the hypotheses from Theorem \ref{aver-theo-2} are fulfilled.
\par First of all, for a given bounded set $D$, it is easy to see that $F^j(r,\theta)$ and $R^j(r,\theta,\eps)$, $j=1,2$, are Lipschitz with respect to the variable $r$ since they are polynomials in $r$; moreover, they are also periodic on $\theta$ with period $2 \pi$, hence the hypothesis \textit{H2} is satisfied. We are left to show that so are the hypotheses \textit{H1} and \textit{H3} --- this will be done by proving the following propositions:
\begin{prop}\label{proph1}
For $|\eps|\neq 0$ sufficiently small, there exists an open bounded set $C$ such that every solution of system \eqref{main-syst} reaches $\sig$ on $\sig^c$.
\end{prop}
\noindent \textbf{Proof of Proposition \ref{proph1}:} Since the original system is autonomous, it's sufficient to analyze under which conditions the set $\sig^c$ exists.
\par First, notice that we can write $\sig=h^{-1}(0)$, where $h(x,y)=y$. Let $X$ and $Y$ denote the smooth pieces of the system, i.e.:
$$X(x,y)=\begin{pmatrix}
    y\\-x-\eps(f(x)\cdot y+g(x))
\end{pmatrix} \text{ and } Y(x,y)=\begin{pmatrix}
    y\\-x-\eps(f(x)\cdot y-g(x))
\end{pmatrix}.$$
\par Let $p\in \sig$ be the point where the solution crosses the discontinuity, then $p=(x,0)$; computing the Lie derivatives on $p$:
\begin{align*}
    & Xh(p)=X(p)\cdot \nabla h(p)=(0,-x+\eps \cdot g(x))\cdot (0,1)=-x+\eps \cdot g(x),\\
    & Yh(p)=Y(p)\cdot \nabla h(p)=(0,-x-\eps \cdot g(x))\cdot (0,1)=-x-\eps \cdot g(x),
\end{align*}
thus $Xh(p)\cdot Yh(p)=x^2-\eps^2 (g(x))^2$. Therefore, for a sufficiently small $|\eps|$ and $x \neq 0$, $p\in \sig ^c$, hence we can find a bounded set such that every solution passing through it reaches $\sig$ at a crossing point. \qed

\begin{prop}\label{proph3}
The coefficients of the polynomials $f$ and $g$ can be chosen in a way that, for every $a\in \R$ with $F_0(a)=0$, there exists a neighborhood $U$ of $a$ such that $F_0(z)\neq 0$ for all $z \in \overline{U}\setminus \{a\}$ and $d_B(F_0,U,0)\neq 0$
\end{prop}
\noindent \textbf{Proof of Proposition \ref{proph3}:} Let $Z=\{r \in \R^+: F_0(r)=0\}$ be the set of positive zeros of the averaged function. The Brouwer degree of a $C^1$-function $f$ with respect to a neighborhood $V$ of zero is given by (see \cite{buica-llibre}):
$$d_B(f,V,0)=\sum_{a\in f^{-1}(0)\cap V} \text{sign(det}Df(a))$$
and, since $F_0$ is a polynomial on $r\in \R$, then
$$d_B(F_0,V,0)=\sum_{a\in Z\cap V} \text{sign}(F_0'(a)).$$
\par Let $a\in Z$, then there exists a neighborhood $U$ such that $F_0(z)\neq 0$ $\forall z \in \overline{U}\setminus \{a\}$ (e.g. consider an open interval with radius half the distance to the next zero of $F_0$); then
$$ d_B(F_0,U,0)= \text{sign}(F_0'(a)).$$
\par Thus what is left to prove is that we can choose the coefficients of $F_0$ such that $F_0'(a)\neq 0$ for every $a\in Z$.
\par Recall that, if $f(x)=\sum_{i=0}^{n}a_ix^i$ and $g(x)=\sum_{j=0}^m b_jx^j$, then
\begin{align*}
    F_0(r)=\frac{1}{2\pi}\left[\sum_{i=0}^{\left[\frac{n}{2}\right]} \pi a_{2i} \alpha_{2i} r^{2i+1} + 2\cdot \sum_{j=0}^{\left[\frac{m}{2}\right]}\frac{-2 b_{2j} r^{2j}}{2j+1} \right],
\end{align*}
or simply
\begin{align*}
    F_0(r)=\sum_{i=0}^{\left[\frac{n}{2}\right]} \Tilde{a}_{2i} r^{2i+1} +  \sum_{j=0}^{\left[\frac{m}{2}\right]} \Tilde{b}_{2j} r^{2j}.
\end{align*}
\par Since $m\geq 1$, we can choose $b_0\neq 0$, hence $\Tilde{b}_0\neq 0$. We will show that this is sufficient to prove that, in order to choose the coefficients of a polynomial with a maximum number of positive roots, we can't have the derivatives vanishing at the roots.
\par Consider $p(x)=c_0+\sum_{j=1}^n c_j x^{l_j}$ a polynomial with $n+1$ terms for which we have chosen the coefficients such that it has $n$ positive roots. Notice that Theorem \ref{theodescartes} guarantees such choice and this number is maximal. Since $p(x)$ is a $ C^{\infty}$-function, then so is $p'(x)$; hence, the Mean Value Theorem implies that between each zero of $p(x)$ there should exist a zero of $p'(x)$. However, $p'(x)=\sum_{j=1}^n l_j\cdot c_j x^{l_j-1}$ is a polynomial with $n$ terms; thus, by Theorem \ref{theodescartes}, it can have at most $n-1$ positive roots, which means that all the positive real roots of $p'(x)$ lie between the positive roots of $p(x)$, therefore $p$ can't have it's derivative vanishing at its positive roots.
\par Since $F_0$ fulfills the above conditions, we can choose $\Tilde{a}_{2i}$ and $\Tilde{b}_{2j}$ so that $F_0(r)$ has $\left[\frac{n}{2}\right]+\left[\frac{m}{2}\right]+1$ roots and for every $a$ with $F_0(a)=0$ there's a neighborhood $U\ni a$ in which $F_0(z)\neq 0$ for all $z\in U$ and
$$d_B(F_0,U,0)=\text{sign}(F_0'(a)) \neq 0.$$
\par Furthermore, since $\Tilde{a}_{2i}=\pi a_{2i}\alpha _{2i}$ with $\alpha_{2i}\neq 0$ and $\Tilde{b}_{2j}=\frac{-4 b_{2j}}{2j+1}$, we can choose $a_{i}$ and $b_j$, $i=1,...,n$ and $j=1,...,m$, such that we obtain the desired coefficients for $F_0$. \qed

\section{Example}

\par Let's find a system like \eqref{main-syst} with $n=4$ and $m=2$ such that the maximum number of predicted limit cycles is achieved, which is, in this case, $4$ limit cycles. In order to do so, consider the following system:
$$\begin{matrix*}[l]
\dot{x}=y,\\
\dot{y}=-x-\eps [(a_0+a_1x+a_2x^2+a_3x^3+a_4x^4)\cdot y +\text{sgn}(y)(b_0+b_1x+b_2x^2)].\end{matrix*}$$
\par Our first step is to make a change of variables to write the original system in polar coordinates, which yields the system:
\begin{align*}
\dot{r}=&-\eps  \sin ( \theta) (r \sin ( \theta) ({a_0}+r \cos ( \theta) ( {a_1}+r \cos ( \theta) ( {a_2}+r \cos ( \theta) ( {a_3}+ {a_4} r \cos ( \theta)))))+\\
&+{\text{sgn}}(r)  {\text{sgn}}(\sin ( \theta)) ( {b_0}+r \cos ( \theta) ( {b_1}+ {b_2} r \cos ( \theta)))),\\
\dot{\theta}=&-\sin ( \theta) (\eps  \cos ( \theta) ( {a_0}+r \cos ( \theta) ( {a_1}+r \cos ( \theta) ( {a_2}+r \cos ( \theta) ( {a_3}+ {a_4} r \cos ( \theta)))))+\sin ( \theta))+\\
&-\frac{\eps   {\text{sgn}}(r) \cos ( \theta)  {\text{sgn}}(\sin ( \theta)) ( {b_0}+r \cos ( \theta) ( {b_1}+ {b_2} r \cos ( \theta)))}{r}-\cos ^2( \theta).\end{align*}
\par Taking $\theta$ as the new independent variable and applying the Taylor expansion series until order 2, we get the following differential equation:
\begin{align*}
    \frac{dr}{d\theta}=&\eps (r \sin ^2( \theta) (\  {a_0}+r \cos ( \theta) (\  {a_1}+r \cos ( \theta) (\  {a_2}+r \cos ( \theta) (\  {a_3}+\  {a_4} r \cos ( \theta)))))+\\
    &+{\text{sgn}}(\sin ( \theta)) (\  { b_0}+r \cos ( \theta) (\  { b_1}+\  { b_2} r \cos ( \theta))))+\mathcal{O}(\eps^2).
\end{align*}
\par Then we can compute the averaged function $F_0(r)=\frac{1}{2\pi}\int_0^{2\pi}F(r,\theta)d\theta$, where $F(r,\theta)$ is the expression multiplying $\eps$ in the above equation:
\begin{align*}
    F_0(r)=\frac{\eps  \left(r \left(3 \pi  \left(8 a_0+2 a_2 r^2+a_4 r^4\right)+32 b_2 r\right)+96 b_0\right)}{48 \pi }.
\end{align*}
\par Having found the averaged function, our next step is to force the existence of four positive roots; i.e. we build a system on the coefficients of $F_0(r)$ by replacing $r$ by four positive values, namely $r=1,2,3 \text{ and } 4$. 
\begin{align}\label{roots-aver-ex}
    \left\{\begin{matrix*}[l]
      3 \pi  (8 a_0+2 a_2+a_4)+32 b_2+96 b_0=0,\\
      2 (3 \pi  (8 a_0+8 a_2+16 a_4)+64 b_2)+96 b_0=0,\\
    3 (3 \pi  (8 a_0+18 a_2+81 a_4)+96 b_2)+96 b_0=0,\\
    4 (3 \pi  (8 a_0+32 a_2+256 a_4)+128 b_2)+96 b_0=0.
    \end{matrix*}\right.
\end{align}
\par Solving the system \eqref{roots-aver-ex} on $a_0$, $a_2$, $a_4$ and $b_0$, we get
\begin{align*}
\left\{a_0 = -\frac{476 b_2}{225 \pi },a_2 = -\frac{52 b_2}{45 \pi },a_4 = \frac{8 b_2}{225 \pi },b_0 = \frac{4 b_2}{15}\right\}.
\end{align*}
\par If we set $b_2=1$ and $a_1=a_3=b_1=0$, the original system will be as following:
\begin{align}\label{ex-4-cycles}
\begin{matrix*}[l]
\dot{x}=y,\\
\dot{y}=-x-\eps  \left(\left(\frac{8 x^4}{225 \pi }-\frac{52 x^2}{45 \pi }-\frac{476}{225 \pi }\right) y+\left(x^2+\frac{4}{15}\right) \text{sgn}(y)\right).  
\end{matrix*}\end{align}
\par
\begin{figure}[h!]
	\centering
    \includegraphics[width=0.55\textwidth]{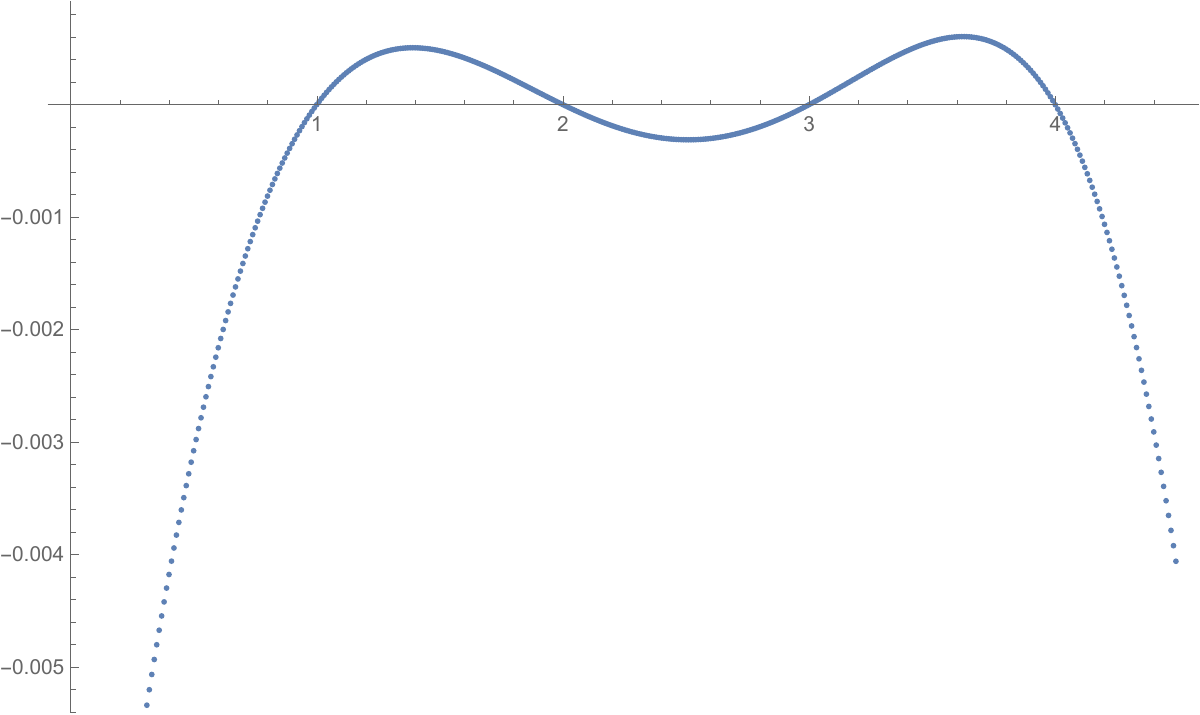}
    \caption{Poincaré Map from the example}
    \label{poincare}
\end{figure}

 If we do the inverse calculations, the averaged system derived from system \eqref{ex-4-cycles} will have the exact 4 roots that we forced in our calculations, which implies by Theorem \ref{aver-theo-2} that, for a sufficiently small $\eps$, this system should have 4 limit cycles. Indeed, using \textit{Wolfram Mathematica} \cite{Mathematica} we are able to detect these limit cycles for $\eps=\frac{1}{100}$ studying the Pincaré map near the points $(1,0)$, $(2,0)$, $(3,0)$ and $(4,0)$. In the Figure \ref{poincare}, the $y-$axis corresponds to the difference $P(x)-x$, where $P(x)$ denotes the image of a point $x$ in the $x-$axis by the Poincaré map. 

\begin{figure}[h!]
  \centering
  \subfloat[Global phase portrait of system \eqref{ex-4-cycles}]{\includegraphics[width=0.4\textwidth]{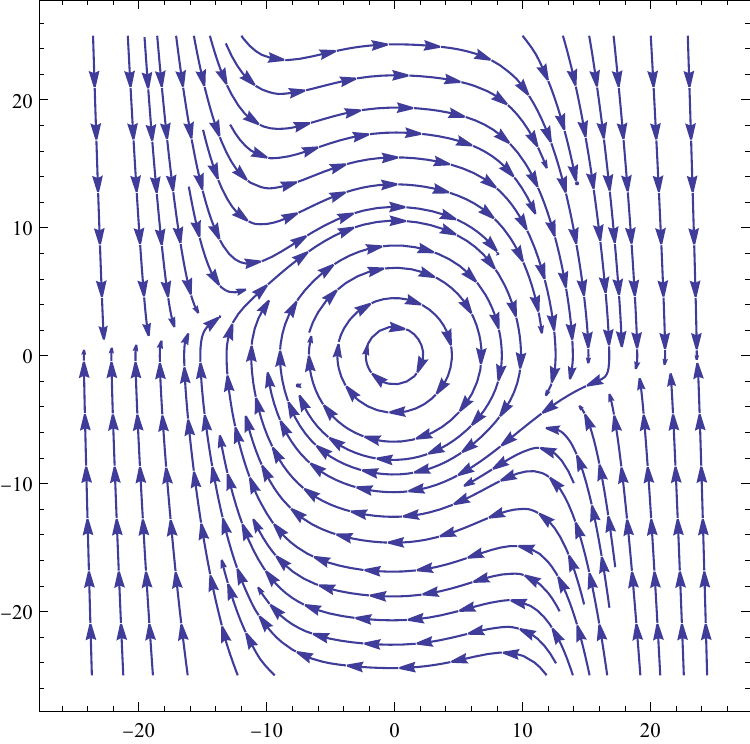}\label{retrato-fase-ex}}
  \hfill
  \subfloat[The green cycles are stable and the red cycles are unstable]{\includegraphics[width=0.4\textwidth]{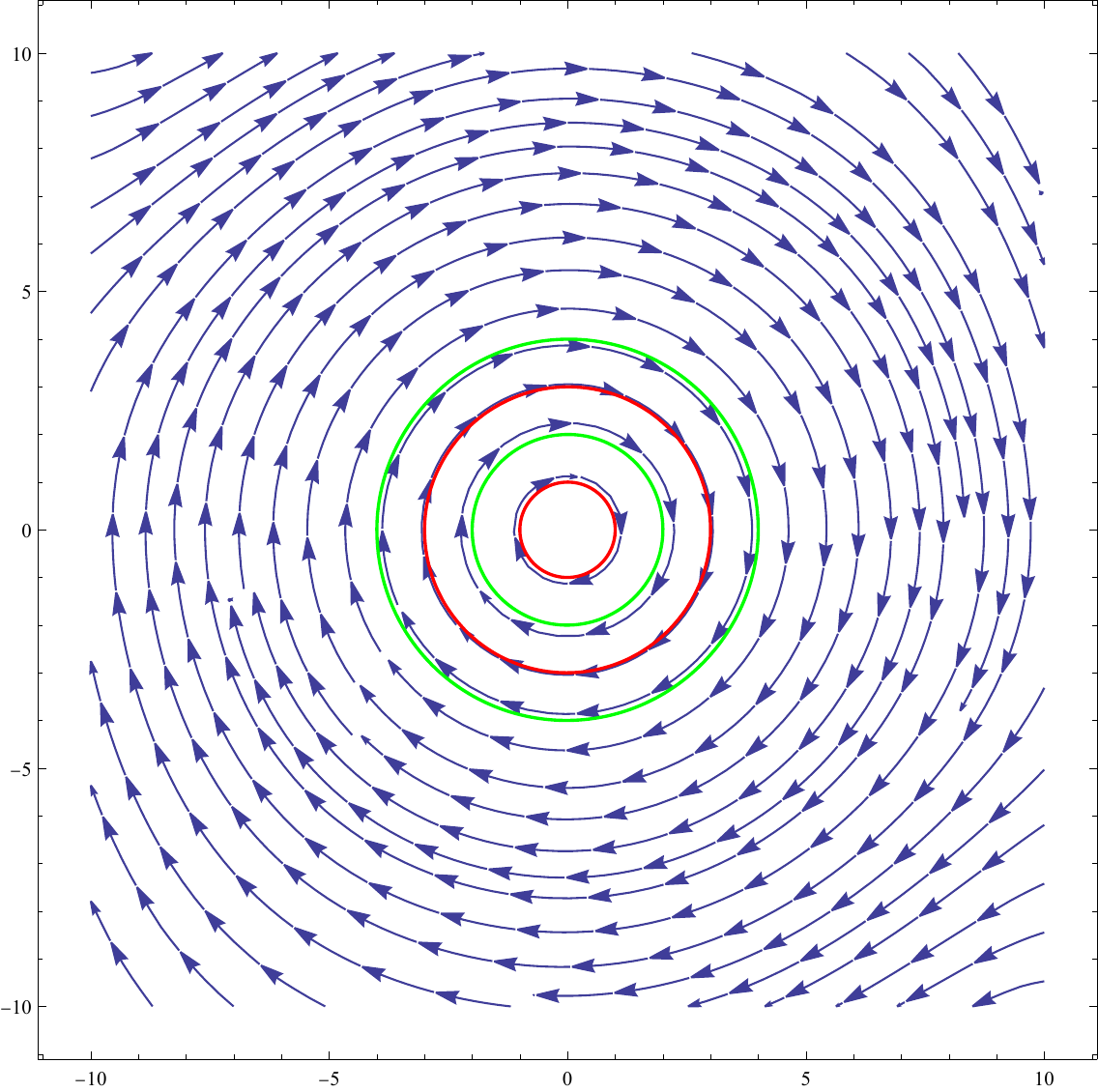}\label{ciclos-ex}}
  \caption{Illustration from the general behavior of system \eqref{ex-4-cycles}}
  \label{last-ex-fig}
\end{figure}

We notice that between the origin and the point $(1,0)$ the solutions are spiraling towards the origin, while between $(1,0)$ and $(2,0)$ they are repelled from the direction of the origin. This behavior will invert itself between $(2,0)$ and $(3,0)$, then again between  $(3,0)$ and $(4,0)$ and then, finally, after $(4,0)$ the solutions are all attracted towards the origin. Therefore we will have an unstable limit cycle near $(1,0)$, a stable limit cycle near $(2,0)$, another unstable one near $(3,0)$ and another stable one near $(4,0)$. Figure \ref{last-ex-fig} illustrates the limit cycles in the phase portrait.

\section*{Acknowledgements}
The São Paulo Research Foundation (FAPESP) partially supports R.M.M. grants 2021/08031-9, 2018/03338-6 and supports T.M.P.A. grant 2022/07654-5. The National Council for Scientific and Technological Development (CNPq) partially supports R. M. M. grants 315925/2021-3 and 434599/2018-2 and partially supports T.M.P.A. grant 132226/2020-0. This study was financed in part by the Coordena\c c\~ao de Aperfei\c coamento de Pessoal de N\'ivel Superior - Brasil (CAPES) - Finance Code 001.

\printbibliography
\end{document}